\documentclass[9pt]{article}
\usepackage[all]{xy}
\usepackage{amsthm,amsfonts,amsmath,amssymb,amscd}
\usepackage{bm}
\usepackage[dvips]{color}
\theoremstyle{definition}
\newtheorem{thm}{Theorem} [section]
\newtheorem{cor}[thm]{Corollary}
\newtheorem{lem}[thm]{Lemma}
\newtheorem{pro}[thm]{Proposition}

\newtheorem{rem}[thm]{Remark}
\newtheorem{defi}[thm]{Definition}
\newtheorem{exm}[thm]{Example}
\author{Hitoshi Yamanaka}
\title{\textbf{Equivariant hyperbolic diffeomorphisms and representation coverings} }
\date{}

\begin{document}
\maketitle
\begin{abstract}
Let $G$ be a compact Lie group and $X$ be a compact smooth $G$-manifold with finitely many $G$-fixed points. 
We show that if $X$ admits a $G$-equivariant hyperbolic diffeomorphism having
a certain convergence property, there exists an
open covering of $X$ indexed by the $G$-fixed points so that each open set is $G$-stable and $G$-equivariantly diffeomorphic to
the tangential $G$-representation at the corresponding $G$-fixed point. We also show that the converse is also true in case of holomorphic torus actions.

\end{abstract}

\section{Introduction}

It is known that there is a certain similarity between algebraic $\mathbb{C}^{*}$-actions on smooth complex projective varieties and 
smooth $\mathbb{R}$-actions on compact Riemannian manifolds arising from the negative gradient flows for Morse functions (see, for example, [CG, 2.4]).
The first similarity is that the both actions induce decompositions of acted spaces via affine spaces, that is, the Bialynicki-Birula decompositions 
and the Morse decompositions.

On the other hand, a theorem of Bialynicki-Birula [BB] states that the existence of an algebraic $\mathbb{C}^{*}$-action on a smooth complex projective
variety implies the existence of an open covering of the variety so that each open set is isomorphic to the affine space having the same dimension and is stable under the
$\mathbb{C}^{*}$-action. So, roughly specking, if a smooth complex projective variety has an algebraic $\mathbb{C}^{*}$-action, it is covered by representation spaces
of $\mathbb{C}^{*}$. 

The aim of the present paper is to give a Morse theoretic counterpart of the above result of Bialynicki-Birula.
More generally, we consider a compact smooth manifold having an invariant hyperbolic diffeomorphism
satisfying a certain convergence condition. Let us explain our result precisely.

Let $G$ be a compact Lie group and $X$ be a compact smooth $G$-manifold with finitely many $G$-fixed points. Let $\varphi :X\rightarrow X$ be a $G$-invariant 
hyperbolic diffeomorphism of $X$. Then we say that the hyperbolic diffeomorphism $\varphi $ satisfying \textbf{convergence condition}
if the following three conditions are satisfied: (1) the $G$-fixed point set coincides with the fixed point set of $\varphi $,
(2) for each $G$-fixed point $p$, the intersection of the corresponding stable and unstable manifolds is the singleton $\{ p\}$ (3)
for any $x\in X$, the sequence $\varphi ^{n}(x)$ converges to a $G$-fixed point when $n$ tends to $\infty$ or $-\infty$.
Note that if $\varphi $ is the hyperbolic diffeomorphism arising from the negative gradient flow of a $G$-invariant Morse function
relative to a $G$-invariant Riemannian metric, it satisfies the convergence condition. In particular, for a compact torus $T$, every compact Hamiltonian 
$T$-space with finitely many $T$-fixed points admits such a hyperbolic diffeomorphism.

We next make the following definition:

\begin{defi}
An open covering $(U_{p}|p\in X^{G})$ of $X$ indexed by $G$-fixed points is called a \textbf{$\bm{G}$-representation covering} of $X$
if each open set $U_{p}$ is stable under the $G$-action and is $G$-equivariantly diffeomorphic to the tangential $G$-representation $T_{p}X$.
\end{defi}

Under the above terminology, we will prove the following theorem:
\begin{thm}
Assume that the compact smooth $G$-manifold $X$ has a $G$-equivariant hyperbolic diffeomorphism satisfying the convergence condition.
Then $X$ admits a $G$-representation covering. 
\end{thm}

Theorem 1.2 implies the following results:

\begin{cor}
If $X$ has a $G$-invariant Morse function, there exists a $G$-representation covering of $X$.
\end{cor}
\begin{cor}
Every compact Hamiltonian $T$-space with finitely many $T$-fixed points admits a $T$-representation covering.
\end{cor}

\bigskip

Concerning Theorem 1.2, a natural question arises: Does the converse of Theorem 1.2 hold ?
The following second theorem gives a partial affirmative answer, that is, the converse is also true in case of holomorphic torus actions:

\begin{thm}
Let $T$ be a compact torus. Assume that a compact complex manifold $X$ with holomorphic $T$-action admits a $T$-representation 
covering as a holomorphic chart. Then it admits a $T$-equivariant hyperbolic diffeomorphism satisfying the convergence condition.
\end{thm}
\bigskip

The contents of the present paper is as follow. In Section 2 we call some basic definitions concerning hyperbolic diffeomorphisms
and hyperbolic dynamical systems. In Section 3, we prove Theorem 1.2. In Section 4, we give a method to construct a hyperbolic
dynamical system from torus actions and invariant almost complex structures. Finally, in Section 5 we prove Theorem 1.5.
\bigskip
\\
{\bf Acknowledgment.} The work was supported by Grant-in-Aid for JPSP Fellows $23\cdot 9693$.

\section{Hyperbolic dynamical systems}

In this section we recall some fundamental definitions concerning hyperbolic dynamical systems and 
hyperbolic diffeomorphisms. For the details, we refer [BH], [PD].
\bigskip

Let $X$ be a compact manifold and $\xi $ be a smooth vector field on $X$. Let $\varphi :\mathbb{R}\times X\rightarrow X$ be the corresponding 
dynamical system. So the curve $\{\varphi (s,x)\}_{s\in \mathbb{R}}$ gives the flow of $\xi $ whose initial value is $x$.

Let $p$ be a fixed point of $\varphi $. Recall that the differential $(d\xi )_{p}:T_{p}X\rightarrow T_{p}X$ is an $\mathbb{R}$-linear transformation defined as follow:
for a local frame $\xi _{1},\dots,\xi _{n}\in \Gamma (TU)\ (n=\dim X)$ on an open neighborhood $U$ of $p$, the representation matrix of $(d\xi )_{p}$ with respect to the ordered base $(\xi _{1})_{p},
\dots,(\xi _{n})_{p}$ is given by $[(\xi _{i})_{p}f_{j}]_{i,j}$ where $f_{i}:U\rightarrow \mathbb{R}$ are $C^{\infty}$-function satisfying
$\xi |_{U}=\sum_{1\leq i\leq n}f_{i}\xi _{i}$. It is straightforward to check that the definition is independent of choice of local frame.

\begin{defi}
\begin{itemize}
\item[(1)] A fixed point $p$ of $\varphi $ is called a \textbf{hyperbolic fixed point} of $\varphi $ if the differential $(d\xi )_{p}$ has
no complex eigenvalues whose absolute values are $1$. 
\item[(2)] We say that $\xi $ is a \textbf{hyperbolic vector field} if all fixed points of $\varphi $ are hyperbolic.
\end{itemize}
\end{defi}

Let $\xi $ be a hyperbolic vector field and $p$ be a hyperbolic fixed point. The \textbf{stable} and \textbf{unstable manifold} of $p$ are
defined by
\[
W^{s}(p):=\bigg\{ x\in X \bigg| \lim\limits_{s\rightarrow \infty}\varphi (s,x)=p \bigg\} 
\]
and
\[
W^{u}(p):=\bigg\{ x \in X \bigg| \lim\limits_{s\rightarrow -\infty}\varphi (s,x)=p \bigg\}
\]
respectively. 

We next explain the notion of hyperbolic diffeomorphisms. Let $\varphi :X\rightarrow X$ be a diffeomorphism of $X$.

\begin{defi}
\begin{itemize}
\item[(1)] A fixed point $p$ of $\varphi $ is called \textbf{hyperbolic fixed point} of $\varphi $ if the differential $(d\varphi )_{p}:T_{p}X\rightarrow T_{p}X$
has no complex eigenvalues whose real parts are zero.
\item[(2)] We say that $\varphi $ is a \textbf{hyperbolic diffeomorphism} of $X$ if all fixed points of $\varphi $ are hyperbolic.
\end{itemize}
\end{defi}

Let $\varphi $ be a hyperbolic diffeomorpshism of $X$ and $p$ be a hyperbolic fixed point. The \textbf{stable} and \textbf{unstable manifold} of $p$ are
defined by
\[
W^{s}(p):=\bigg\{ x\in X \bigg| \lim\limits_{n\rightarrow \infty}\varphi ^{n}(x)=p \bigg\} 
\]
and
\[
W^{u}(p):=\bigg\{ x \in X \bigg| \lim\limits_{n\rightarrow -\infty}\varphi ^{n}(x)=p \bigg\}
\]
respectively. The \textbf{stable manifold theorem} states that they are injectively immersed submanifold of $X$.

We summarize the facts which will be needed in the later sections:
\begin{lem} Let $\xi $ be a vector field on $X$ and $\varphi :\mathbb{R}\times X\rightarrow X$ be the corresponding
dynamical system.
\begin{itemize}  
\item[(1)] $\xi $ is hyperbolic if and only if the diffeomorphism $\varphi_{s}:X\rightarrow X$ defined by $\varphi _{s}(x):=\varphi (s,x)$ is
hyperbolic for all $s\in \mathbb{R}\setminus \{0\}$.

\item[(2)] If the limits $\lim\limits_{s\rightarrow \infty}\varphi (s,x)$ and $\lim\limits_{s\rightarrow -\infty}\varphi (s,x)$ exist for all $x\in X$,
then the limits $\lim\limits_{n\rightarrow \infty}\varphi_{s} ^{n}(x)=p$ and  $\lim\limits_{n\rightarrow -\infty}\varphi_{s} ^{n}(x)=p $ exist for all $s\in \mathbb{R},x\in X$. Also, the corresponding stable and unstable manifolds coincide.
\end{itemize}
\begin{proof}
For (1), we refer [BH, p.113, Lemma 4.19]. (2) is clear.
\end{proof}
\label{dyn}
\end{lem}

\section{Representation covering}

Let $G$ be a compact Lie group and $X$ be a compact smooth $G$-manifold with finitely many $G$-fixed points. 
We introduce our main subject in this paper, that is, the notion of $G$-representation coverings of $X$:

\begin{defi}
We say that a smooth $G$-manifold $X$ has a \textbf{$\bm{G}$-representation covering} if there exists an open covering
$(U_{p}|p\in X^{G})$ of $X$ indexed by $G$-fixed points such that
\begin{itemize}
\item[(1)] each $U_{p}$ is $G$-stable,
\item[(2)] $U_{p}\cong T_{p}X$ as $G$-manifolds.
\end{itemize}

\end{defi}

We give three examples of smooth $G$-manifolds admitting representation coverings.
\begin{exm}
Let 
\[
X=S^{2n}=\{\ (z_{1},\dots,z_{n},s)\in \mathbb{C}^{n}\times \mathbb{R}\ |\ |z_{1}|^{2}+\dots+|z_{n}|^{2}+s^{2}=1\ \}
\]
be the $2n$-dimensional
sphere. We define an action of $T=(S^{1})^{n}$ on $S^{2n}$ by
\[
(t_{1},\dots,t_{n})\cdot (z_{1},\dots,z_{n},s)=(t_{1}z_{1},\dots,t_{n}z_{n},s).
\]
Then one can construct a $T$-representation covering of the pair $(S^{2n},T)$ as follow:

Let $S^{2n}=U_{0}\cup U_{\infty}, U_{0}=\{s\not=1\},U_{\infty}=\{s\not= -1\}$ be the standard $T$-invariant open covering. 

Then the maps 
\[
\varphi _{0}:U_{0}\rightarrow T_{(0,\cdots,0,1)}S^{2n}
\]
\[\varphi _{\infty}:U_{\infty}\rightarrow T_{(0,\cdots,0,-1)}S^{2n}
\]
given by
\[
\varphi _{0}(z_{1},\dots,z_{n},s)=\big(\frac{\overline{z_{1}}}{1-s},\dots,\frac{\overline{z_{n}}}{1-s}\big),
\]
\[
\varphi _{\infty}(z_{1},\dots,z_{n},s)=\big(\frac{z_{1}}{1+s},\dots,\frac{z_{n}}{1+s}\big)
\]
are $T$-equivariant. Note that the former is the composition of the stereographic projection and the componentwise complex conjugation
$(z_{1},\cdots,z_{n})\mapsto (\overline{z_{1}},\dots,\overline{z_{n}})$.
\end{exm}

\begin{exm}
Let 
\[
X=\mathbb{C}P^{n}=\{\ [z_{0}:z_{1}:\dots:z_{n}]\ |\ (z_{0},z_{1},\dots,z_{n})\in \mathbb{C}^{n+1}\setminus \{ (0,0,\dots,0)\}\ \}
\]
be the complex projective space of complex dimension $n$. We define an action of $T=(S^{1})^{n}$ by
\[
(t_{1},\dots,t_{n})\cdot [z_{0}:z_{1}:\dots:z_{n}]=[z_{0}:t_{1}z_{1}:\dots:t_{n}z_{n}].
\]

Then one can construct a $T$-representation covering of $\mathbb{C}P^{n}$ as follow:

Let $U_{i}$ be an open set of $\mathbb{C}P^{n}$ defined by
\[
U_{i}:=\{\ [z_{0}:z_{1}:\dots:z_{n}]\in \mathbb{C}P^{n}\ |\ z_{i}\not=0\ \}\ \ \ \ \ \ \ \ \ (0\leq i\leq n).
\]

Then the open covering $(U_{i}|0\leq i\leq n)$ gives a $T$-representation covering. The same construction is possible for arbitrary toric manifold, i.e., complete non-singular toric variety with the Hausdorff topology (see [F]).
\end{exm}
\begin{exm} 
Let $G$ be a connected reductive Lie group over $\mathbb{C}$ and $B$ be a Borel subgroup of $G$. We take $X$ and $T$ as the flag manifold $G/B$ and the compact real form of the 
maxiaml torus of $G$ contained in $B$, respectively.

Then one can construct a $T$-representation covering of the pair $(G/B,T)$ as follow: Let $X(w_{0})$ be the Bruhat cell corresponding to the longest element $w_{0}$
of $W$. We set $U_{w}:=ww_{0}^{-1}X(w_{0})\ (w\in W)$. Then [S, p.149, 8.5.1. Proposition (ii)] and [S, p.152, 8.5.10. Excercise (1)] imply that the open sets $\{U_{w}|w\in W\}$ of $X$ gives an open covering of the flag manifold $X$ and each $U_{w}$ is $T$-equivariantly diffeomorphic to 
\[
\bigoplus _{\alpha \in \Delta ^{+}}\mathfrak{g}_{-w.\alpha }.
\]
Here, $\Delta ^{+}$ is the set of positive roots corresponding to the Borel subgroup $B$ and $\mathfrak{g}_{-w\alpha }$ is the root space for the root $-w.\alpha $.

\end{exm}

To state our result, we make the following definition:
\begin{defi}
We say that a $G$-invariant hyperbolic diffeomorphism 
\[
\varphi :X\longrightarrow  X
\]
has \textbf{convergence property} if it satisfies the following
three properties:
\begin{itemize}
\item[(1)] $X^{G}\subset \text{Fix}(\varphi )$,
\item[(2)] $W^{u}(p)\cap W^{s}(p)=\{ p\}$ for all $p\in X^{G}$,
\item[(3)] the limits $\lim\limits_{n\rightarrow \infty}\varphi^{n}(x)$ and $\lim\limits_{n\rightarrow -\infty}\varphi^{n}(x)$ converge
to some $G$-fixed points of $X$ for all $x\in X$. 
\end{itemize}
\end{defi}

The following is our main theorem in this section:

\begin{thm}
If $X$ admits a $G$-equivariant hyperbolic diffeomorphism having convergence property, then there exists a $G$-representation covering of
$X$.
\end{thm}

As a corollary one finds that the existence of a $G$-representation covering gives an obstruction for the existence of
$G$-invariant Morse functions:

\begin{cor}
If $X$ admits a $G$-invariant Morse function, there exists a $G$-representation covering of $X$.
\begin{proof}
Take a $G$-invariant Riemannian metric on $X$. Let $\widetilde{\varphi} :\mathbb{R}\times X\rightarrow X$ be the hyperbolic dynamical system associated to the 
negative gradient vector field of the $G$-invariant Morse function
relative to the $G$-invariant Riemannian metric. Then the hyperbolic diffeomorphism $\varphi :X\rightarrow X$ defined by $\varphi (x):=\widetilde{\varphi }(1,x)$
gives a $G$-equivariant hyperbolic diffeomorphism having convergence property.
\end{proof}
\end{cor}

Corollary 3.7 implies the following. Let $T$ be a compact torus.
\begin{cor}
Every compact Hamiltonian $T$-manifold $X$ with finitely many $T$-fixed points admits a $T$-representation covering.
\begin{proof}
Let $\Phi :X\rightarrow (\text{Lie}(T))^{*}$ be the moment map. Then for a generic $a \in \text{Lie}(T)$, the contraction $\Phi _{a}:X\rightarrow \mathbb{R}$,
$\Phi _{a}(x):=(\Phi (x))(a)$ give a $T$-invariant Morse function as well known (see for example [MS2]).
\end{proof}
\end{cor}

\begin{rem}
It is known a similarity between algebraic $\mathbb{C}^{*}$-actions on smooth complex projective varieties and the $\mathbb{R}$-actions arising from negative gradient flows associated to Morse functions on compact
manifolds (for example, see [CG, 2.4]). 

A theorem of Bialynicki-Birula [BB] states that every smooth complex projective variety with an algebraic $\mathbb{C}^{*}$-action
has a Zariski open covering so that each Zariski open set is $\mathbb{C}^{*}$-stable and is isomorphic to an algebraic $\mathbb{C}^{*}$-representation.

From this point of view, Theorem 3.6 (or Corollary 3.7) can be thought as a Morse theoretic counterpart of the above theorem of Bialynicki-Birula.

\end{rem}

The rest of this section is devoted to the proof of Theorem 3.6. We begin by the following lemma:

\begin{lem}
Let $E$ be a Banach space, $U$ be an open neighborhood of $0\in E$ and $f:U\rightarrow E$ be a $C^{\infty}$-map such that
$f(0)=0,(df)_{0}=\text{id}_{E}$. If a positive real number $r>0$ satisfies that by measuring the operator norm, $|\text{id}_{E}-(df)_{x}|<1/2$ for any
$x\in B_{2r}(0)$, then the restriction of $f$ induces a $C^{\infty}$-diffeomorphism
\[
f^{-1}(B_{r/2}(0))\cap B_{r}(0)\longrightarrow B_{r/2}(0).
\]

\begin{proof}
This follows from the proof of the inverse mapping theorem given in [L].
\end{proof}
\label{inv}
\end{lem}

Let $N$ be a $G$-manifold (not necessary compact), $p$ be a $G$-fixed point of $N$ and $\varphi :N\rightarrow N$ be a $G$-equivariant hyperbolic diffeomorphism of $N$ having
the following properties:
\begin{itemize}
\item $N^{G}=\text{Fix}(\varphi )=\{ p\}$
\item $p$ is the global attractor, i.e., $N=W^{s}(p)$.
\end{itemize}

We fix a $G$-invariant Riemannian metric of $N$ and denote by $\|\cdot \|$ the
corresponding fiberwise norm of $TN$. For a positive real number $r>0$, we denote by $B_{r}(0)$ the open ball in $T_{p}N$ which is centered at $0$ and have radius $r$ with respect to 
the norm $\|\cdot \|$. We note that since the norm $\|\cdot \|$ is $G$-invariant, each $B_{r}(0)$ is stable under the $G$-action on $T_{p}N$. The set $\{B_{r}(0)|r>0\}$ forms
a fundamental neighborhood system of $0\in T_{p}N$.

By focusing on the exponential map with respect to the $G$-invariant Riemannian metric, one can find a positive real number $R>0$,
a $T$-invariant open neighborhood $U$ of $p$ and a $G$-equivariant diffeomorphism
\[
\psi :U\longrightarrow B_{R}(0)\ (\subset T_{p}N).
\]

We define a $G$-equivariant diffeomorphism 
\[
g:\psi (U\cap \varphi ^{-1}(U))\longrightarrow \psi (U\cap \varphi (U))
\]
by
\[
g:=\psi |_{U\cap \varphi (U)}\circ \varphi |_{U\cap \varphi ^{-1}(U)}\circ \psi ^{-1}|_{\psi (U\cap \varphi ^{-1}(U))}.
\]

In the rest of this section we set $E:=T_{p}N$ for simplicity. We note that the vector space $E$ can be viewed as a $G$-representation, i.e.,
the tangential representation at $p$. Let us define a $C^{\infty}$-map $u:E\rightarrow E$ by
$u(\xi )=\rho (\|\xi \|)\xi $ where $\rho$ is a $C^{\infty} $-function $\rho :\mathbb{R}\rightarrow \mathbb{R}$ depicted 
by the following picture:
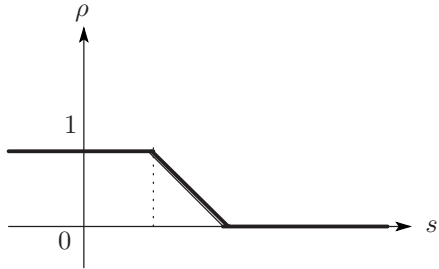
\[
\unitlength 0.1in
\begin{picture}( 21.2800, 14.0200)( 14.1000,-24.6000)
% VECTOR 2 0 3 0 Black White
% 2 1410 2244 3509 2244
% 
{\color[named]{Black}{%
\special{pn 8}%
\special{pa 1410 2244}%
\special{pa 3510 2244}%
\special{fp}%
\special{sh 1}%
\special{pa 3510 2244}%
\special{pa 3442 2224}%
\special{pa 3456 2244}%
\special{pa 3442 2264}%
\special{pa 3510 2244}%
\special{fp}%
}}%
% LINE 0 0 3 0 Black White
% 2 1410 1850 2167 1850
% 
{\color[named]{Black}{%
\special{pn 20}%
\special{pa 1410 1850}%
\special{pa 2168 1850}%
\special{fp}%
}}%
% LINE 0 0 3 0 Black White
% 2 2531 2244 3394 2244
% 
{\color[named]{Black}{%
\special{pn 20}%
\special{pa 2532 2244}%
\special{pa 3394 2244}%
\special{fp}%
}}%
% STR 2 0 3 0 Black White
% 4 1793 1021 1793 1123 5 0 0 0
% $\rho$
\put(17.9300,-11.2300){\makebox(0,0){$\rho$}}%
% STR 2 0 3 0 Black White
% 4 3623 2141 3623 2244 5 0 0 0
% $s$
\put(36.2300,-22.4400){\makebox(0,0){$s$}}%
% STR 2 0 3 0 Black White
% 4 1670 2257 1670 2360 2 0 0 0
% $0$
\put(16.7000,-23.6000){\makebox(0,0)[lb]{$0$}}%
% STR 2 0 3 0 Black White
% 4 1736 1606 1736 1709 5 0 0 0
% $1$
\put(17.3600,-17.0900){\makebox(0,0){$1$}}%
% LINE 2 2 3 0 Black White
% 2 2167 1832 2167 2244
% 
{\color[named]{Black}{%
\special{pn 8}%
\special{pa 2168 1832}%
\special{pa 2168 2244}%
\special{dt 0.045}%
}}%
% VECTOR 2 0 3 0 Black White
% 2 1803 2460 1803 1206
% 
{\color[named]{Black}{%
\special{pn 8}%
\special{pa 1804 2460}%
\special{pa 1804 1206}%
\special{fp}%
\special{sh 1}%
\special{pa 1804 1206}%
\special{pa 1784 1274}%
\special{pa 1804 1260}%
\special{pa 1824 1274}%
\special{pa 1804 1206}%
\special{fp}%
}}%
% LINE 0 0 3 0 Black White
% 4 2160 1850 2557 2241 2557 2241 2557 2241
% 
{\color[named]{Black}{%
\special{pn 20}%
\special{pa 2160 1850}%
\special{pa 2558 2242}%
\special{fp}%
\special{pa 2558 2242}%
\special{pa 2558 2242}%
\special{fp}%
}}%
% LINE 2 0 3 0 Black White
% 2 2530 2240 2150 1860
% 
{\color[named]{Black}{%
\special{pn 8}%
\special{pa 2530 2240}%
\special{pa 2150 1860}%
\special{fp}%
}}%
\end{picture}%
\]

We note that since the fiberwise norm $\|\cdot \|$ is $G$-invariant and since $G$ acts on $E$ as a representation, the $C^{\infty}$-map $u$ is $G$-equivariant.

Let $A$ be the differential $(dg)_{0}:E\rightarrow E$ of $g$ at $0$ and we set 
\[
\Phi :=A^{-1}\circ g:\psi (U\cap \varphi ^{-1}(U))\longrightarrow E.
\]

The both of $A$ and $\Phi $ are $G$-equivariant as easily understood (note that since $G$ acts on $E$ as a $G$-representation, $E\cong T_{0}E$ as $G$-representations). In the rest of this section, we set $V:=\psi (U\cap \varphi ^{-1}(U))$
for simplicity.

What we first want to show is that the germ of the map $g:\psi (U\cap \varphi ^{-1}(U))\rightarrow E$ at $0$ coincides with
the germ of a $G$-equivariant diffeomorphism $h:E\rightarrow E$ at $0$, that is, the maps coincide on an open neighborhood of $0$.

For the proof of Theorem 3.5, we need to go back to the proof of the so-called stable manifold theorem. 
We should remark that our arguments in Lemma 3.11 and Lemma 3.12 are greatly influenced on the paper of Abbondandolo-Mayer [AM].

Under the above notations, we next show the following lemma:

\begin{lem}
There exist positive real numbers $r>0,\delta >0$ and $C^{\infty}$-maps $X:[0,1]\times B_{r/2}(0)\rightarrow E,\ Y:[0,1]\times B_{\delta }(0)\rightarrow E$
satisfying the following conditions:
\begin{itemize}
\item[(1)] $Y([0,1]\times B_{\delta }(0))\subset B_{r/2}(0)$,
\item[(2)] $Y(0,\xi )=\xi$ and $\displaystyle \frac{\partial Y}{\partial t}(t,\xi )=X(t,G(t,\xi ))$ for all $(t,\xi )\in [0,1]\times B_{\delta }(0)$,
\item[(3)] $Y(1,\xi )=\Phi (\xi )$ for all $\xi \in B_{\delta }(0)$,
\item[(4)] $X$ and $Y$ are $G$-equivariant with respect to the second components,
\item[(5)] $X(t,0)=0, (D_{2}X_{t})_{0}=0$ for each $t\in [0,1]$. Here $X_{t}$ is given by $X_{t}(\xi )=X(t,\xi )$.
\end{itemize}
\begin{proof}
We follow up [AM, Lemma 4.8]. We define a $C^{\infty}$-map $Y:[0,1]\times V\rightarrow E$ by $Y(t,\xi )=t\Phi (\xi )+(1-t)\xi $ and set $Y_{t}(\xi ):=Y(t,\xi )$.
Note that since the map $\Phi $ is $G$-equivariant and $G$ acts on $E$ as a real $G$-representation, the map $Y_{t}$ is $G$-equivariant
for each $t\in [0,1]$.

Let us consider a continuous map $\beta :[0,1]\times V\rightarrow \mathbb{R}$ given by $\beta  (t,\xi ):=|\text{id}_{E}-(dY_{t})_{\xi }|$.
Since $(dY_{t})_{0}=\text{id}_{E}$ for any $t\in [0,1]$, one finds that the set $\beta ^{-1}(\{s\in \mathbb{R}|s<1/2\})$ is an open neighborhood
of $[0,1]\times \{0\}$ in $[0,1]\times V$. In particular, for each $t\in [0,1]$, we can choose an open neighborhood $U_{t}$ of $t$ in $[0,1]$ and
a positive real number $r_{t}>0$ so that $U_{t}\times B_{2r_{t}}(0)\subset \beta ^{-1}(\{s\in \mathbb{R}|s<1/2\})$. Thanks to the compactness
of the closed interval $[0,1]$, there exist finitely many real numbers $t_{1},\cdots,t_{n}\in [0,1]$ so that $[0,1]$ is covered by
$\{U_{t_{i}}|1\leq i\leq n\}$. We define a positive real number $r>0$ by $r:=\min_{1\leq i\leq n}r_{t_{i}}$. Then one can easily deduce that 
$[0,1]\times B_{2r}(0)\subset \beta ^{-1}(\{s\in \mathbb{R}|s<1/2\})$.

Therefore we have that $| \text{id}_{E}-(dY_{t})_{\xi }|<1/2$ for all $(t,\xi )\in [0,1]\times B_{2r}(0)$. Thus by Lemma \ref{inv}, we can conclude that
there exists a positive real number $r>0$ which is independent of the parameter $t\in [0,1]$ so that for each $t\in [0,1]$, the restriction of $Y_{t}$ induces
a $C^{\infty}$-diffeomorphism 
\[
Y_{t}^{-1}(B_{r/2}(0))\cap B_{r}(0)\longrightarrow B_{r/2}(0)
\]
We can choose the real number $r$ so that $B_{r}(0)\subset V$ by taking $r$ smaller if necessary.

Denote by $\overline{Y}_{t}$ the above $C^{\infty}$-diffeomorphism and define a $C^{\infty}$-map $H:[0,1]\times B_{r/2}(0)\rightarrow E$
by $H(t,\xi ):=\overline{Y}_{t}(\xi )$. It is easy to check that $H$ is $T$-equivariant with respect to the second component.

We next claim that 
\[
\bigcap _{t\in [0,1]}Y_{t}^{-1}(B_{r/2}(0))=Y_{0}^{-1}(B_{r/2}(0))\cap Y_{1}^{-1}(B_{r/2}(0)).
\]

Clearly the LHS is contained in the RHS. If $\xi $ belongs to the RHS, one has $\xi=Y_{0}(\xi )\in B_{r/2}(0)$ and
$\Phi (\xi )=Y_{1}(\xi )\in B_{r/2}(0)$. Thus by the convexity of the open ball $B_{r/2}(0)$, we have that $Y_{t}(\xi )=t\Phi (\xi )+(1-t)\xi \in B_{r/2}(0)$
for all $t\in [0,1]$, as desired.

As a consequence, there exists a positive real number $\delta >0$ so that
\[
B_{\delta }(0)\subset \bigcap _{t\in [0,1]}Y_{t}^{-1}(B_{r/2}(0))\cap B_{r}(0)\cap V.
\]

Then each $\overline{Y}_{t}$ induces an injection $B_{\delta }(0)\rightarrow B_{r/2}(0)$. In particular, for each $(t,\xi )\in [0,1]\times B_{\delta }(0)$
we have $Y(t,\xi )\in B_{r/2}(0)$ and have $H(t,Y(t,\xi ))=\xi $. 

Since we have chosen $r>0$ so that $B_{r}(0)\subset V$, the image of $H:[0,1]\times B_{\delta }(0)\rightarrow E$ is contained in $V$. This allows us to
define a $C^{\infty}$-map $X:[0,1]\times B_{r/2}(0)\rightarrow E$ by $X(t,\xi )=\Phi (H(t,\xi ))-H(t,\xi )$. It is clear that $X$ is $G$-equivariant with
respect to the second component.

In the above arguments, we have constructed two $C^{\infty}$-maps
\[
X:[0,1]\times B_{r/2}(0)\longrightarrow E,\ \ \ \ \ \ \ \ \ \ Y:[0,1]\times B_{\delta }(0)\longrightarrow E.
\]

It is straightforward to check that the positive real numbers $r,\delta$ and the above two maps $X,Y$ satisfy our
desired properties.
\label{XY}
\end{proof}
\end{lem}

\begin{lem}
There exists a positive real number $\epsilon >0$ and a $G$-equivariant $C^{\infty}$-diffeomorphism $h:E\rightarrow E$ satisfying
the following conditions:
\begin{itemize}
\item[(1)] $B_{\epsilon }(0)\subset V\cap B_{R}(0)$,
\item[(2)] $h|_{B_{\epsilon  }(0)}=g|_{B_{\epsilon   }(0)}$, 
\item[(3)] for each $x\in N$, there exists a non-negative integer $n(x)$ so that $h^{n(x)}(x)\in B_{\epsilon   }(0)$.
\end{itemize}
\begin{proof}
We follow up [AM, Proposition 4.7]. Let $X:[0,1]\times B_{r/2}(0)\rightarrow E$ and $Y:[0,1]\times B_{\delta }(0)\rightarrow E$
as in Lemma \ref{XY}.  

Since the function $\rho $ is identity on an open neighborhood of $0\in \mathbb{R}$ and has compact support, there exist
two positive real numbers $0<s_{0}<r_{0}$ so that $\varphi |_{B_{s_{0}}}=\text{id}_{B_{s_{0}}(0)}$ and $\varphi (E)\subset B_{r_{0}}(0)$.
By the condition (5) in Lemma 3.11 and the Tayler's formula, the argument presented in Lemma 3.11 shows that there exist positive real numbers $\epsilon _{1},r_{1}>0$ so that $e^{\epsilon _{1}\frac{r_{0}}{s_{0}}}\| A\|<1$ and $|X(t,\xi )|\leq \epsilon _{1}|\xi |$ for all $(t,\xi )\in [0,1]\times B_{r_{1}}(0)$. Let us define $\widetilde{u}:E\rightarrow E$ and $\widetilde{X}:[0,1]\times E\rightarrow E$ by
\[
\widetilde{u}(\xi )=\frac{r_{1}}{r_{0}}u \big(\frac{r_{0}}{r_{1}}\xi \big), \ \ \ \ \ \ \ \widetilde{X}(t,\xi )=X(t,\widetilde{u }(\xi )).
\]
Note that $\widetilde{u}|_{B_{s}(0)}=\text{id}_{B_{s}(0)}$ for $s:=(r_{1}s_{0})/r_{0}$. Also note that $\widetilde{u}$ is $G$-equivariant and $\widetilde{X}$ is $G$-equivariant with respect to the second component.

We claim that $|\widetilde{X}(t,\xi )|\leq \epsilon _{1}(r_{1}/s_{0})|\xi| $ for all $(t,\xi )\in [0,1]\times E$. If
if $\xi\in B_{s_{0}}(0) $, the estimate folds since $X(t,\xi )=\widetilde{X}(t,\xi) $. If $\xi\not\in B_{s_{0}}(0)$, we have that $|\widetilde{X}(t,\xi )|\leq \epsilon _{1}|\widetilde{u }|\leq \epsilon _{1}r_{1}\leq \epsilon _{1}(r_{1}/s)|\xi |=\epsilon _{1}(r_{0}/s_{0})|\xi |$.
As a consequence, for each $c>0$, $\widetilde{X}(t,\xi )$ is bounded velocity on $[0,1]\times B_{c}(0)$ in sense of [H, p.178].

Thus, [H, p.179, 1.1.Theorem] implies that there exists a unique diffeotopy $\widetilde{Y}:[0,1]\times E\rightarrow E$. So we have that $(\partial Y /\partial t)(t,\xi )=\widetilde{X}(t,\widetilde{G}(t,\xi ))$ and $\widetilde{Y}(0,\xi )=\xi $. Thanks to the $G$-equivariancy of $\widetilde{X}$, the map $g^{-1}\cdot \widetilde{Y}(t,g\cdot \xi )$ is also the solution of this equation for each $g\in G$.
Thus one finds that $\widetilde{Y}$ is $G$-equivariant with respect to the second component. We also note that the estimate for $\widetilde{X}$ implies that $|\widetilde{Y}(t,\xi )|
\leq e^{\epsilon _{1}(r_{0}/s_{0})t}|\xi |$ for all $(t,\xi )\in [0,1]\times E$.
Moreover, since $X$ and $\widetilde{X}$ coincide locally, the argument presented in Lemma 3.11 shows that there exists a positive real number $\epsilon >0$ so that $Y$ and $\widetilde{Y}$ coincide on $[0,1]\times B_{\epsilon }(0)$, and $B_{\epsilon}(0)\subset V\cap B_{R}(0)\cap B_{\delta }(0)$. 

Let us define a diffeomorphism $\widetilde{\Phi }:E\rightarrow E$ by $\widetilde{\Phi }(\xi ):=\widetilde{Y}(1,\xi )$. Then we claim that a map $h:E\rightarrow E$ defined by
$h:=T\circ \widetilde{\Phi }$ is the desired diffeomorphism. Note that $h$ is clearly $G$-equivariant.

The condition (1) is clear. The condition (2) folds since $Y(1,\xi )=\Phi (\xi )$ for all $\xi \in B_{\delta }(0)$ and $Y$ and $\widetilde{Y}$ coincide locally.
Finally, since 
\[
|h(\xi )|\leq \| A\| |\widetilde{\Phi }(\xi )|\leq \|A \| e^{\epsilon _{1}(r_{0}/s_{0})}|\xi |
\]
and $\|A \| e^{\epsilon _{1}(r_{0}/s_{0})}<1$, the condition (3) holds.

The proof is now complete.

\end{proof}
\label{loc}
\end{lem}

\begin{pro}
There exists a $G$-equivariant $C^{\infty}$-diffeomorphism $T_{p}N\rightarrow N$.
\begin{proof}
This follows from a standard argument and Lemma 3.12. We describe the proof to check the $G$-equivariancy.

Take a $G$-equivariant $C^{\infty}$-diffeomorphism $h:E\rightarrow E$ as in Lemma \ref{loc}.
Let us define a map $\phi :N\rightarrow E$ as follow. For a point $x\in N$, one can find a non-negative integer
$n$ so that $\varphi^{n}(x)\in \psi ^{-1}(B_{\epsilon }(0))$. By using the integer $n$, we define $\phi (x):=h^{-n}(\psi (\varphi ^{n}(x)))$. We claim that
the value $\phi (x)$ is independent of the integer $n$ whenever $\varphi^{n}(x)\in \psi ^{-1}(B_{\epsilon }(0))$. To prove the claim, we first note that
$\psi (\varphi ^{n}(x))\in B_{\epsilon }(0)$.  Since we have chosen
the differmorphism $h$ so that $h|_{B_{\epsilon }(0)}=g|_{B_{\epsilon }(0)}$, one obtains that $h(\psi (\varphi ^{n}(x)))=g(\psi (\varphi ^{n}(x)))$.
Henceforth we have that
\begin{align*}
h^{-n}(\psi (\varphi ^{n}(x))) &= h^{-(n+1)}(h(\psi (\varphi ^{n}(x)))) \\
                               &= h^{-(n+1)}(g(\psi (\varphi ^{n}(x)))) \\
                               &= h^{-(n+1)}((\psi\circ \varphi\circ \psi ^{-1})(\psi (\varphi ^{n}(x)))) \\
                               &= h^{-(n+1)}(\psi(\varphi ^{n+1}(x))).
\end{align*}
This completes the proof of the well-definedness of $\phi $. Note that $\phi $ is clearly $G$-equivariant.

We next construct the inverse of $\phi $. For $\xi \in E$, one can find a non-negative integer $n$ so that $h^{n}(\xi )\in B_{\epsilon }(0)$.
Then we define $\overline{\phi}:E\rightarrow N$ by $\overline{\phi }(\xi ):=\varphi ^{-n}(\psi^{-1}( h^{n}(\xi )))$.
The same argument shows that the value $\overline{\phi }(\xi )$ is independent of the integer $n$ whenever $h^{n}(\xi )\in B_{\epsilon }(0)$.

Since $\overline{\phi }$ clearly gives the inverse of $\phi $, the proof is now complete.

\end{proof}
\label{vec}
\end{pro}

We go to the next step. Let $W$ be a finite dimensional real $G$-representation and $\pi:E\rightarrow W$ be a smooth $G$-equivariant vector bundle over $W$. Let $\nabla:\Gamma (TW)\times \Gamma (E)\rightarrow \Gamma (E) $ be a $G$-equivariant 
covariant derivative, that is, a covariant derivative having the property that 
\[
(\nabla_{g\cdot \xi }(g\cdot s))(g\cdot x)=g\cdot ((\nabla_{\xi }s)(x))
\]
for all $g\in G,\xi \in \Gamma (TW),s\in \Gamma (E),x\in N$. Note that such a covariant derivative exists since for a fixed covariant derivative $\nabla $,
one can construct a $G$-equivariant covariant derivative $\widetilde{\nabla}$ by averaging  as follow:
\[
(\widetilde{\nabla}_{\xi }s)(x):=\int_{g\in G}[g^{-1}\cdot \nabla _{g\cdot \xi }(g\cdot s)](x)dg.
\]
Here $dg$ is the Haar measure of $G$ normalized so that $\int_{g\in G}dg=1$.
 
For a point $x\in W$, we define a curve $c_{x}:[0,1]\rightarrow W$ by $c_{x}(t):=(1-t)x$. We denote by $P_{c_{x}}:\pi^{-1}(x)\rightarrow \pi^{-1}(0)$ the parallel
transport associated to the curve $c_{x}$ and the covariant derivative $\nabla $. Let $\Gamma (c_{x})$ be the set of smooth sections of $E$ along the curve $c_{x}$.
Note that an element $g\in G$ induces a map $g_{*} :\Gamma (c_{x})\rightarrow \Gamma (c_{g\cdot x})$ defined by $(g_{*}s)(t):=g\cdot (s(t))$.

Recall from elementary differential geometry that there exists an $\mathbb{R}$-linear transformation
\[
\frac{D}{dt}:\Gamma (c_{x})\rightarrow \Gamma (c_{x})
\]
which is characterized by the following two properties:
\begin{itemize}
\item[(1)] $\displaystyle \frac{D}{dt}(fs)=\frac{df}{dt}s+f\frac{Ds}{dt}$ \ \ \ \ \ \ \ \ \ \ \ $(f\in C^{\infty}([0,1]),s\in \Gamma (c_{x}))$,
\item[(2)] if $s\in \Gamma (c_{x})$ is given by $s(t)=\widetilde{s}(c_{x}(t))\ (t\in (a,b))$ for some $\widetilde{s}\in \Gamma (E)$ and 
two real numbers $0\leq a<b\leq 1$, then we have
\[
\frac{Ds}{dt}(t)=\nabla _{c_{x}'(t)}\widetilde{s}\ \ \ \ \ \ \ \ \ (t\in (a,b)).
\]
\end{itemize}
\begin{lem}
The diagram 
\[
\begin{xy}
\xymatrix{ \ar@{}[rd]| {}
\Gamma (c_{x}) \ar[r]^{\frac{D}{dt}} \ar[d]_{g_{*}} & \Gamma (c_{x}) \ar[d]^{g_{*}}\\
\Gamma (c_{g\cdot x}) \ar[r]_{\frac{D}{dt}} & \Gamma(c_{g\cdot x}) }
\end{xy}
\]
commutes for all $g\in G,x\in W$.
\label{comm}
\begin{proof}
Take $s\in \Gamma (c_{x})$ and $t_{0}\in [0,1]$. Then by considering a local frame of $E$ around $c_{x}(t_{0})$,
one can take a global section $\widetilde{s}\in \Gamma (E)$ and two real numbers $0\leq a<b\leq 1$ so that $s(t)=\widetilde{s}(c_{x}(t))$ for all $t\in (a,b)$.

Then one has that
\begin{align*}
(g_{*}\frac{Ds}{dt})(t_{0})  &= g\cdot (\nabla _{c_{x}'(t_{0})}\widetilde{s})                        \\
                             &= \nabla _{c_{g\cdot x}'(t_{0})}(g\cdot \widetilde{s})\ \ \ \ \ \ \ \ \ \ \ \ \ \ \ \ \ \ \ \ (\text{as\ $\nabla$\ is\ $G$-equivariant}).            
\end{align*}
On the other hand, since
\[
(g\cdot \widetilde{s})(c_{g\cdot x}(t))=g\cdot (\widetilde{s}(g^{-1}\cdot (c_{g\cdot x}(t))))=g\cdot (\widetilde{s}(c_{x}(t)))=g\cdot (s(t))=(g_{*}s)(t),
\]
we have that
\[
\nabla _{c_{g\cdot x}'(t_{0})}(g\cdot \widetilde{s})=\frac{D(g_{*}s)}{dt}(t_{0}).
\]
This completes the proof.

\end{proof}
\end{lem}

\begin{cor}
We have that $P_{c_{g\cdot x}}(g\cdot v)=g\cdot (P_{c_{x}}(v))$ for all $g\in G,x\in W,v\in \pi^{-1}(x)$.
\begin{proof}
Let $s\in \Gamma (c_{x})$ be the parallel section along to the curve $c_{x}$ (so one has $Ds/dt=0$).

By the definition of the parallel transport, we have that $s(0)=v$ and $P_{c_{x}}(v)=s(1)$. 

Then the curve $g_{*}s\in \Gamma (c_{g\cdot x})$ satisfies that
\[
(g_{*}s)(0)=g\cdot (s(0))=g\cdot v
\]
and 
\[
\frac{D(g_{*}s)}{dt}=g_{*}\frac{Ds}{dt}=0.
\]
by Lemma \ref{comm}.

Thus we have that
\[
P_{c_{g\cdot x}}(g\cdot v)=(g_{*}s)(1)=g\cdot (s(1))=g\cdot (P_{c_{x}}(v))
\]
as desired.
\end{proof}
\label{para}
\end{cor}

\begin{pro}
$E$ is isomorphic to $W\times \pi^{-1}(0)$ as a $G$-equivariant smooth vector bundle over $W$. Here the $G$-action on $E\times \pi^{-1}(0)$
is given by the diagonal one.
\begin{proof}
We define a map $\Psi :E\rightarrow W\times \pi^{-1}(0)$ by $\Psi (u):=(\pi(u),P_{c_{\pi(u)}}(u))$. Then the map $\Psi $ gives an isomorphism of smooth
vector bundles and is $G$-equivariant by Lemma \ref{para}.

\end{proof}
\label{equibdl}
\end{pro}

We are now in the position to prove Theorem 3.6.

\begin{proof}[Proof of Theorem 3.6]
By the assumption (2), $W^{s}(p)$ and $W^{u}(p)$ are embedded $G$-stable submanifold of $X$ (see, for example, [BH, p.115, Lemma 4.20]). Thus one can apply $G$-equivariant tubular neighborhood theorem [K, p.178, Theorem 4.8] to $W^{s}(p)$, and has a $G$-stable open neighborhood $U_{p}$ of $W^{s}(p)$ in $X$
which is $G$-equivariantly diffeomorphic to the normal bundle 
\[
\nu _{p}:=TX|_{W^{s}(p)}/TW^{s}(p)
\]
of $W^{s}(p)$ in $X$.

On the other hand, we have the following series of $G$-equivariant diffeomorphisms:
\begin{align*}
\nu _{p}  &\cong        TX|_{W^{s}(p)}/TW^{s}(p) \\
          &\cong        W^{s}(p)\times (T_{p}X/T_{p}W^{s}(p))\ \ \ \ \ \ \ \ \ \ \ \ (\text{as Proposition \ref{vec} and \ref{equibdl}}) \\
          &\cong        W^{s}(p)\times T_{p}W^{u}(p)       \ \ \ \ \ \ \ \ \ \ \ \ \ \ \ \ \ \ \ \ \ (\text{as}\ T_{p}X=T_{p}W^{s}(p)\oplus T_{p}W^{u}(p)) \\
          &\cong        T_{p}W^{s}(p)\times T_{p}W^{u}(p) \ \ \ \ \ \ \ \ \ \ \ \ \ \ \ \ \ \      (\text{as Proposition \ref{vec}}) \\
          &\cong        T_{p}X.
\end{align*}

Since $X$ is decomposed into stable manifolds $W^{s}(p)$ by the assumption (3), the family $(U_{p}|p\in X^{G})$ of $G$-stable open sets gives
a desired $G$-representation covering of $X$.
\end{proof}

\section{Topological generators and hyperbolicity}

In the previous section, we have shown that if a smooth $G$-manifold admits a certain equivariant hyperbolic diffeomorphism,
it also admits a $G$-representation covering. 

Then the following natural question arise: ``Does the converse of Theorem 3.14 hold ?".
Unfortunately, the author do not have any satisfactory answer to the question.

The aim of this and the next section is to give a partial affirmative answer for the above question. Roughly specking, we will show that the converse is also true
in case of holomorphic torus actions.
\bigskip

In this section we give a method to construct hyperbolic dynamical systems from torus actions with invariant almost complex structures.

Let $T=(S^{1})^{r}$ be the compact torus of rank $r$ and $X$ be a smooth $2n$-dimensional $T$-manifold with an almost complex structure $J$ compatible with the $T$-action (the term``compatible" means that for each $t\in T,x\in X$, the induced map $t_{*}:T_{x}X\rightarrow T_{t\cdot x}X$ is $\mathbb{C}$-linear with respect to
$J_{x}$ and $J_{t\cdot x}$).

To construct a hyperbolic dynamical system from the torus action, we focus on a topological generator $t_{0}$ of $T$ (by the definition, the cyclic group $\{ t_{0}^{k}|k\in \mathbb{Z}\}$ generated by $t_{0}$ is a dense subset of the compact torus $T$). We also take an element $a_{0}$ of the Lie algebra $\mathfrak{t}$ of $T$ so that $\exp(a_{0})=t_{0}$.
Let us denote by $\xi _{0}$ the fundamental vector field associated to $a_{0}$. By the definition we have
\[
(\xi _{0})_{x}\phi =\frac{d}{ds}\phi (\exp (-sa_{0})\cdot x)\bigg|_{s=0}
\]
for all smooth function $\phi:X\rightarrow \mathbb{R}$ and $x\in X$.

We set $\xi _{0}^{J}:=-J\xi _{0}$, and denote by $\text{zero}(\xi_0)$ and $\text{zero}(\xi _{0}^{J})$ the set of zero points of $\xi _{0}$ and $\xi _{0}^{J}$ respectively.
\begin{lem}
We have the following:
\[
\text{zero}(\xi_0^{J})=\text{zero}(\xi _{0})=X^{T}.
\]
\begin{proof}
The first equality and the inclusion $\text{zero}(\xi _{0})\supset X^{T}$ are clear.
To see the inverse inclusion $\text{zero}(\xi _{0})\subset  X^{T}$, let $p\in \text{zero}\ (\xi _{0})$. Since $t_{0}$ is a topological generator of $T$,
it is enough to show that $t_{0}\cdot p=p$.

Assume that $t_{0}\cdot p\not= p$. Then there exists a smooth function $\phi :X\rightarrow \mathbb{R}$ which separates $t_{0}\cdot p$ and $p$, that is, $\phi (t_{0}\cdot p)\not=\phi (p)$.

For a real number $s_{0}$, we define a smooth function $\phi _{0}:\mathbb{R}\rightarrow \mathbb{R}$ by
\[
\phi _{0}(s)=\phi \Big(\exp \big(-sa_{0}\big)\cdot p\Big).
\]
Then we have
\begin{align*}
\frac{d}{ds}\phi _{0}(s)\bigg| _{s=s_{0}} &= \frac{d}{ds}\phi (\exp (-sa_{0})\cdot p)\bigg| _{s=s_{0}} \\
                                          &= \frac{d}{ds}\phi (\exp (-(s+s_{0})a_{0})\cdot p)\bigg| _{s=0} \\
                                          &= \frac{d}{ds}(\phi\circ \exp (-s_{0}a_0))(\exp (-sa_{0})\cdot p)\bigg| _{s=0} \\
                                          &= (\xi _{0})_{p}(\phi\circ \exp(-s_{0}a_0)) \\
                                          &= 0.
\end{align*}
Thus the function $\phi _{0}$ is a constant function. In particular, by taking the values at $s=0$ and $s=-1$, we have
$\phi (t_{0}\cdot p)=\phi (p)$. This is a contradiction.

\end{proof}
\end{lem}
\begin{rem}
The following seems to be well-known:
\begin{center}
if $H^{\text{odd}}(X)=\{0\}$, we have that $X^{T}\not= \emptyset$.
\end{center}
In fact, one can give a short proof of this fact using our vector field $\xi _{0}$:
Assume that $X^{T}=\emptyset$. Then $\text{zero}(\xi _{0})=\emptyset$ by Lemma 4.1. Henceforth [MS1, PROPERTY 9.7]
implies that the Euler class of $X$ vanishes. So the Euler characteristic of $X$ also vanishes. This
is a contradiction.
\end{rem}

Let us denote by $\varphi ^{J}:\mathbb{R}\times X\rightarrow X$ the dynamical system corresponding
to $\xi _{0}^{J}$.

We take a Riemannian metric $g$ of $X$ which is $T$-invariant and $J$-Hermitian (note that for a $T$-invariant metric $g_{0}$, 
the metric
\[
g(u,v):=\frac{g_{0}(u,v)+g_{0}(Ju,Jv)}{2}\ \ \ \ (u,v\in T_{x}X,x\in X)
\]
has the required property since we assumed that $J$ is compatible with the $T$-action). For a $T$-fixed point $p$ of $X$, let $\exp_{p}:T_{p}X\rightarrow X$ be the exponential map
associated to the Levi-Cibita connection of $g$. This is $T$-invariant and the differential at the origin of $T_{p}X$ is the identity map
under the natural identification $T_0(T_{p}X)=T_{p}X$. Thus there exists a positive real number $\delta $ so that the restriction
\[
\exp_{p}|_{B_{\delta }(0)}:B_{\delta }(0)\rightarrow \exp_{p}(B_{\delta }(0))
\]
 gives a $T$-equivariant diffeomorphism (here $B_{\delta }(0)$ is the
open ball whose center is 0 and radius is $\delta $ in $T_{p}X$ with respect to the metric $g$). We set $U_{p}:=\exp_{p}(B_{\delta }(0))$.

Moreover, by composing the diffeomorphism and the $T$-equivariant diffeomorphism
\[
T_{p}X\rightarrow B_{\delta }(0),\ \ \ u\mapsto \frac{\delta u}{\sqrt{1+\| u \|^{2}}},
\]
we obtain a $T$-equivariant diffeomorphism $\varphi _{p}:T_{p}X\rightarrow U_{p}$. 

Let 
\[
T_{p}X=\bigoplus_{i=1}^{n}(T_{p}X)^{\alpha _{p,i}}
\]
be the irreducible decomposition of the $T$-representation space $T_{p}X$. Here $\alpha _{p,i}:T\rightarrow S^{1}$ is
a weight of the $T$-representation and 
\[
(T_{p}X)^{\alpha _{p,i}}=\{ u\in T_{p}X|\ t\cdot u=\alpha _{p,i}(t)u\ (t\in T)\ \}
\] 
is the corresponding weight space.

We fix a $\mathbb{C}$-base $u_{p,1},\cdots,u_{p,n}$ of $T_{p}X$ so that $u_{p,i}\in (T_{p}X)^{\alpha _{p,i}}$ for all $i$. Then we have a local coordinate system $(U_{p};x_{1},y_{1},\cdots,x_{n},y_{n})$ around $p$ defined by
\[
\varphi _{p}^{-1}(x)=\sum_{i=1}^{n}(x_{i}(x)+\sqrt{-1}y_{i}(x))u_{p,i}.
\]

\begin{lem}
Assume that $0<|X^{T}|<\infty$. Then $\varphi ^{J}:\mathbb{R}\times X\rightarrow X$ is a hyperbolic dynamical system on $X$. 
\begin{proof}
We set
\[
\xi_{0}\big|_{U_{p}}=\sum_{i=1}^{n}f_{i}\frac{\partial }{\partial x^{i}}+\sum_{i=1}^{n}g_{i}\frac{\partial }{\partial y^{i}} \ \ \ \ (f_{i},g_{i}\in C^{\infty }(U_{p})).
\]
So one has
\[ 
\xi_{0}^{J}\big|_{U_{p}}=-\sum_{i=1}^{n}f_{i}J\frac{\partial }{\partial x^{i}}-\sum_{i=1}^{n}g_{i}J\frac{\partial }{\partial y^{i}}.
\]

By the definition of the coordinate functions $x_{1},y_{1},\dots,x_{n},y_{n}$ and $T$-equvaliancy of $\varphi _{p}:T_{p}X\rightarrow U_{p}$, we have 
\begin{align*}
f_{i}(x) &= (\xi_{0})_{x}x_{i} \\
         &=\frac{d}{ds}x_{i}(\exp (-sa_{0})\cdot x)\bigg|_{s=0}\\
         &=\frac{d}{ds}\bigg( \text{$\mathbb{R}$-coefficient of}\ u_{p,i}\ \text{in}\ \exp (-sa_{0})\cdot \bigg(\sum_{i=1}^{n}(x_{i}(x)+\sqrt{-1}y_{i}(x))u_{p,i}\bigg)\bigg)\bigg|_{s=0} \\
         &=\frac{d}{ds}\Big(x_{i}(x)\text{Re}\ \alpha _{p,i}\big(\exp (-sa_{0})\big)-y_{i}(x)\text{Im}\ \alpha _{p,i}\big(\exp (-sa_{0})\big)\Big) \bigg|_{s=0} \\
         &=-x_{i}(x)\text{Re}\ \alpha' _{p,i}\big(a_{0}\big)+y_{i}(x)\text{Im}\ \alpha _{p,i}'\big(a_{0}\big) \\
         &=y_{i}(x)\text{Im}\ \alpha _{p,i}'\big(a_{0}\big)
\end{align*}
for all $x\in U_{p}$.
A similar calculation shows that
\[
g_{i}(x)=-x_{i}(x)\text{Im}\ \alpha _{p,i}'\big(a_{0}\big).
\]
Thus we obtain that
\[
\xi_{0}^{J}\big|_{U_{p}}=-\sum_{i=1}^{n}y_{i}\text{Im}\ \alpha _{p,i}'\big(a_{0}\big)J\frac{\partial }
{\partial x^{i}}+\sum_{i=1}^{n}x_{i}\text{Im}\ \alpha _{p,i}'\big(a_{0}\big)J\frac{\partial }{\partial y^{i}}.
\]

By using the local frame 
\[
\Bigg\{ J\frac{\partial }{\partial x^{1}},J\frac{\partial }{\partial y^{1}},\cdots,J\frac{\partial }{\partial x^{n}},J\frac{\partial }{\partial y^{n}} \Bigg\},
\]
to compute the eigenvalues of $(d\xi_{0}^{J} )_{p}:T_{p}X\rightarrow T_{p}X$, one finds that the representation matrix is given by
\[
\text{diag} (a_{1},a_{1},\cdots,a_{n},a_{n}) \ \ \ \ (a_{i}=\text{Im}\ \alpha _{p,i}'\big(a_{0}\big)).
\]

Let us assume that $\varphi ^{J}$ is not hyperbolic. Then by Lemma 2.3 the $\mathbb{R}$-linear transformation $(d\xi _{0}^{J})_{p}$ has a complex eigenvalue whose real part is zero.
Hence we have $a_{i}=0$ for some $i\in \{ 1,\dots,n\}$. Since $\alpha _{p,i}'$ is $\sqrt{-1}\mathbb{R}$-valued, we also obtain that $\alpha _{p,i}'(a_{0})=0$.
This implies that $\alpha _{p,i}(t_{0})=\alpha _{p,i}(\exp(a_{0}))=\exp(\alpha _{p,i}'(a_{0}))=1$. So $\alpha _{p,i}$ is identically equal to $1$ since
$t_{0}$ is a topological generator of $T$. This implies that $(T_{p}X)^{\alpha _{p,i}}$ is the trivial irreducible $T$-representation.
Since $\exp_{p}:B_{\delta}(0)\rightarrow U_{p}$ is a $T$-equivariant diffeomorphism, we have $X^{T}\supset \exp_{p}(B_{\delta}(0)\cap (T_{p}X)^{\alpha _{p,i}})$.
This is a contradiction since $\dim X^{T}=0$ and $\dim B_{\delta}(0)\cap (T_{p}X)^{\alpha _{p,i}}=2$.

\end{proof}
\end{lem}

\begin{lem}
$\varphi _{s}^{J}(t\cdot x)=t\cdot \varphi _{s}^{J}(x)$ for all $s\in \mathbb{R},t\in T,x\in X$. 
\begin{proof}
Since the almost complex structure $J$ is compatible with the $T$-action, we have
\[
\frac{d}{ds}\varphi _{s}^{J}(t\cdot x)=t_{*}\frac{d}{ds}\varphi _{s}^{J}(x)=-t_{*}(J\xi _{0})_{\varphi _{s}^{J}(x)}
=-J_{t\cdot \varphi _{s}^{J}(x)}(\xi _{0})_{t\cdot \varphi _{s}^{J}(x)}=(\xi _{0}^{J})_{t\cdot \varphi _{s}^{J}(x)}.
\]

This completes the proof.
\end{proof}
\end{lem}

\section{Existence of equivariant hyperbolic diffeomorphisms}

Let $X$ be a compact complex $T$-manifold with finitely many $T$-fixed points admitting a $T$-representation covering $(U_{p}|p\in X^{T})$.
We say that the $T$-representation covering is \textbf{holomorphic} if $U_{p}\cong T_{p}X$ as complex $T$-manifolds. 

In the rest of this section, we assume that the $T$-representation covering is holomorphic.

Let $p$ be a $T$-fixed point of $X$ and $(U_{p};x^{1},y^{1},\cdots,x^{n},y^{n})$ be the coordinate neighborhood induced from the corresponding
complex $T$-representation as in Section 4.

Using the coordinate on $U_{p}$, we define a Riemannian metric $g^{(p)}$ on $U_{p}$ as follow:
\[
g^{(p)}\Big(\frac{\partial }{\partial x^{i}},\frac{\partial }{\partial y^{j}}\Big)=0\ \ \ \ \ \ (1\leq i,j\leq n),
\]
\[
g^{(p)}\Big(\frac{\partial }{\partial x^{i}},\frac{\partial }{\partial x^{j}}\Big)=\delta _{ij}e^{-|x^{i}|},\ \ \ \ \ \ 
g^{(p)}\Big(\frac{\partial }{\partial y^{i}},\frac{\partial }{\partial y^{j}}\Big)=\delta _{ij}e^{-|y^{i}|}\ \ \ \ \ (1\leq i,j\leq n).
\]
Here the symbol $\delta _{ij}$ stands for the Kronecker's delta.

Let $(\rho^{(p)}|p\in X^{T})$ be a partition of unity associated to the open covering $(U_{p}|p\in X^{T})$. Then we define
a Riemannian metric $g$ on $X$ by 
\[
g:=\sum_{p\in X^{T}}\rho^{(p)}g^{(p)}.
\]

We set 
\[
\|u\|^{(p)}:=\sqrt{g^{(p)}(u,u)},\ \ \ \ \ \ \ \|v\|:=\sqrt{g(v,v)}
\]
for $p\in X^{T},u\in TU_{p},v\in TX$.
\begin{pro}
$\|\xi _{0}\|_{\varphi _{s}^{J}(x)}$ converges to 0 when $s$ tends to $\infty$ for any $x\in X$.

\begin{proof}

One can take a $T$-fixed point $p$ so that $x$ is contained in $U_{p}$.

Then the same calculation in the proof of Lemma 4.3 shows that
\[
\xi _{0}^{J}\big|_{U_{p}}=\sum_{i=1}^{n}a_{i}x^{i}\frac{\partial }{\partial x^{i}}+\sum_{i=1}^{n}a_{i}y^{i}\frac{\partial }{\partial y^{i}}
\]
in the coordinate system $(U_{p};x^{1},y^{1},\cdots,x^{n},y^{n})$. 

Henceforth the flow $\varphi _{s}^{J}(x)$ is expressed as follow:
\[
\varphi _{s}(x)=(p_{1}e^{-a_{1}s},q_{1}e^{-a_{1}s},\cdots,p_{n}e^{-a_{n}s},q_{n}e^{-a_{n}s}).
\]
Here $(p_{1},q_{1},\cdots,p_{n},q_{n})\in \mathbb{R}^{2n}$ is the coordinate of the point $x$ in the coordinate neighborhood $(U_{p},;x^{1},y^{1},\cdots,x^{n},y^{n})$.

As a consequence, we have that
\[
\Big(\| \xi _{0}^{J}\|_{\varphi _{s}^{J}(x)}^{(p)}\Big)^{2}=\sum_{i=1}^{n}a_{i}^{2}p_{i}^{2}e^{-2a_{i}s-|p_{i}|e^{-a_{i}s}}+\sum_{i=1}^{n}a_{i}^{2}q_{i}^{2}e^{-2a_{i}s-|q_{i}|e^{-a_{i}s}}.
\]

We claim that when $s\rightarrow \infty$, the term $a_{i}^{2}p_{i}^{2}e^{-2a_{i}s-|p_{i}|e^{-a_{i}s}}$ converges to $0$.
This is because that 
\begin{itemize}
\item if $p_{i}=0$, the claim is obvious.
\item if $p_{i}\not=0$ and $a_{i}>0$, the claim is also obvious.
\item if $p_{i}\not=0$ and $a_{i}<0$, the claim follows from the fact that
\[
-2a_{i}s-|p_{i}|e^{-a_{i}s}\rightarrow -\infty
\]
when $s\rightarrow \infty$.
\end{itemize}

The same consideration shows that the term $a_{i}^{2}q_{i}^{2}e^{-2a_{i}s-|q_{i}|e^{-a_{i}s}}$ also converges to $0$ when $s\rightarrow \infty$.
Thus we have that $\| \xi _{0}^{J}\|_{\varphi _{s}^{J}(x)}^{(p)}\rightarrow 0$ when $s\rightarrow \infty$.

Since $\rho ^{(p)}$ is bounded, we also have $\| \xi _{0}^{J}\|_{\varphi _{s}^{J}(x)}\rightarrow 0$ when $s\rightarrow \infty$.

\end{proof}

\end{pro}

The following is our main theorem in this section:

\begin{thm}
There exists a $T$-invariant hyperbolic diffeomorphism $\varphi :X\rightarrow X$ satisfying the convergence property.
\begin{proof}
The local expression of the vector field $\xi _{0}^{J}$ on $U_{p}$ presented in Lemma 5.1 shows that each open set $U_{p}$ is stable under the 
induced $\mathbb{R}$-action. The expression also implies that $W^{u}(p)\cap W^{s}(p)=\{p\}$ for any $p\in X^{T}$.
So, it is enough to prove that $\varphi ^{J}$ satisfies the condition (3). 

We first consider the case that $s\rightarrow \infty$. Let $x$ be a point of $X$.
Since $X$ is sequentially compact, there exists a real sequence $\{s_{n}\}_{n=1}^{\infty}$ having the following two properties:
\begin{itemize}
\item $\{s_{n}\}_{n=1}^{\infty}$ is strictly increasing and is divergent to $\infty$.
\item $\varphi _{s_{n}}^{J}(x)$ converges to a point $p$ of $X$ when $n$ tends to $\infty$.
\end{itemize}

By Lemma 4.2 we have
\[
\|\xi _{0}^{J}\|_{p}=\|\xi _{0}^{J}\|_{\lim\limits_{n\rightarrow \infty}\varphi _{s_{n}}^{J}(x)}
=\lim\limits_{n\rightarrow \infty}\|\xi _{0}^{J}\|_{\varphi _{s_{n}}^{J}(x)}=0.
\]

Henceforth $p$ is a $T$-fixed point of $X$ by Lemma 4.1. We show that the flow $\varphi_{s}^{J}(x)$ converges to the $T$-fixed point $p$ when $s$ tends to $\infty$.
We proceed proof by contradiction.

Assume that the flow $\varphi_{s}^{J}(x)$ does not converge to the $T$-fixed point $p$. Since $X^{T}$ is finite, one can take a compact neighborhood $K$ of $p$ so that $K\cap X^{T}=\{ p\}$. 
Then there exists an open set $U$ of $X$ having the following two properties:
\begin{itemize}
\item $p\in U\subset K$.
\item There exists a strictly increasing sequence $\{s_{n}'\}_{n=1}^{\infty}$ so that $s_{n}'$ is divergent to $\infty$ and $\varphi _{s_{n}'}^{J}(x)$ is in $K\setminus U$
for all $n\geq 1$.
\end{itemize}
 
Since $X$ is sequentially compact, by taking a subsequence if necessary, one may assume that the sequence $\{ \varphi _{s_{n}'}^{J}(x)\}_{n=1}^{\infty}$
converges to a point $p'$ of $X$. Then the argument used in the proof that $p$ is a $T$-fixed point also implies that $p'$ is a $T$-fixed point of $X$. However, since $\varphi _{s_{n}'}^{J}(x)$ is in 
the closed set $K\setminus U$ for all $n\geq 1$, the $T$-fixed point $p'$ must be contained in $K$ and differs from $p$. This contradicts to the fact that
$K\cap X^{T}=\{ p\}$. This completes the proof that the limit $\lim\limits_{s\rightarrow \infty}\varphi _{s}^{J}(x)$ converges to some $T$-fixed point of $X$.

The convergence of $\lim\limits_{s\rightarrow -\infty}\varphi _{s}^{J}(x)$ follows from the above case by changing the element $a_{0}$ to $-a_{0}$.

The proof is now complete.
\end{proof}
\end{thm}

\end{document}